\documentclass[12pt,a4paper]{amsart}
\usepackage{amssymb}
\usepackage{amsfonts,amsmath,amscd}
\usepackage{latexsym}
\usepackage{euscript}
\usepackage[T2A]{fontenc} \usepackage[cp1251]{inputenc}
\usepackage[russian,english]{babel}
\def\e{\varepsilon}
\def\R{{\mathbb R}}
\def\C{{\mathbb C}}
\def\N{{\mathbb N}}
\def\Re{{\textsf{Re}}}
\def\Im{{\textsf{Im}}}
\def\d{\delta}
\def\a{\alpha}
\def\p{\varphi}

\DeclareMathOperator{\supp}{supp}
\DeclareMathOperator{\dist}{dist}
\DeclareMathOperator{\diam}{diam}
\newtheorem{Th}{Theorem}
\newtheorem{Def}{Definition}

\begin{document}

\title{Uniqueness theorems for almost periodic objects}

\author{S.Yu.~Favorov, O.I.Udodova}

\address{Sergii Favorov,
\newline\hphantom{iii}  Karazin's Kharkiv National University
\newline\hphantom{iii} Svobody sq., 4,
\newline\hphantom{iii} 61022, Kharkiv, Ukraine}
\email{sfavorov@gmail.com}

\address{Olga Udodova,
\newline\hphantom{iii}  ???
\newline\hphantom{iii} ?????,
\newline\hphantom{iii} ??????}
\email{udodova-o@ukr.net}

\maketitle {\small
\begin{quote}
\noindent{\bf Abstract.}
Uniqueness theorems are considered for various types of almost periodic objects: functions,  measures, distributions, multisets, holomorphic and meromorphic functions.
\medskip


\end{quote}
}

\bigskip

It easily follows from the definition of almost periodic functions that if the values of two such functions converge at infinity, then these almost periodic functions coincide.
This effect also manifested itself in \cite{KS} for  zeros of holomorphic almost periodic functions, and then in \cite{F3} and \cite{F4} for Fourier quasicrystals and some classes of transformable measures on LCA-groups.

In this note, we discuss this effect in detail, show how it can be strengthened, what form it takes for other almost periodic objects - almost periodic distributions, almost periodic measures, almost periodic multisets, $a$ -points of holomorphic and meromorphic almost periodic functions.

\medskip

\section{Almost periodic functions}\label{S1}

\medskip

We start with the simplest almost periodic object ~ - uniformly almost periodic functions on a finite-dimensional space and on tube sets. The definitions introduced in this section will also be used in subsequent sections.

 Let $B_\C(z^0,R)$ be the open ball $\{z\in\C:\,|z-z^0|<R\}$ in the space $\C^d$, and  $B_\R(x^0,R)$ be the open ball $\{x\in\R:\,|x-x^0|<R\}$ in the space $\R^d$. The tube set $T_K\subset\C^d$ means the set of the form
$$
T_K=\left\{{z = x + iy\in\C^d :x\in\R^d,\;y\in K} \right\},
$$
where $K$ is a compact subset of $\R^d$. Clearly, $\R^d=T_{\{0\}}$. Then $T_\Omega$ means the domain
$$
T_\Omega   = \left\{ {z = x + iy :x \in \R^d ,\;y \in \Omega } \right\}
$$
where $\Omega$ is a domain in $\R^d$, maybe $\Omega=\R^d$. The set $E$ is relatively dense in  $\R^d$, if there exists $R < \infty$ such that each ball $B_\R(x,R)$ intersects with $E$.
By $\# A$ we denote the number of elements of the finite set $A$.

\begin{Def}
A continuous complex-valued function $f(z)$ on a tube set $T_K$ is called almost periodic if for every $\e>0$ the set of its $\e$-almost periods
$$
E_{\e,K}=E_{\e,K}(f)=\{\tau\in\R^d : \sup_{z\in T_K}|f(z+\tau)-f(z)|<\e\}
$$
is relatively dense in $\R^d$.
\end{Def}

It easily follows from this definition that almost periodic functions on $T_K$ are bounded. Less obvious is the following statement:

\begin{Th}\label{T1}( \cite{R})
A continuous function $f(z)$ on $T_K$ is almost periodic iff for any sequence $\{{x_n}\}\subset\R^d$  there is a subsequence $\{{x_{n'}}\}$
such that the functions $f_{n'}(z)=f(z+x_{n'})$ form the fundamental sequence with respect to the uniform convergence on $T_K$.
\end{Th}

\begin{Def}\hspace{-0,5em}
A function $f(z)$ on a tube domain $T_{\Omega}$ is called almost periodic if for every compact set $K \subset \Omega$  its restriction to $T_K$ is almost periodic.
\end{Def}

\begin{Th}\label{T2}(\cite{R})
A continuous function $f(z)$ is almost periodic on a tube domain $T_{\Omega}$ iff for any sequence $\{ {x_n } \}\subset\R^d$  there is a subsequence $\{ {x_{n'} }\}$
such that the functions $f_{n'}(z)=f(z+x_{n'})$ form the fundamental sequence with respect to the uniform convergence on $T_K$ for every $K \subset \Omega$.
 \end{Th}

{\bf Remark.} All these definitions and theorems carry over practically unchanged to the case of mappings $ F: T_K \to \C^N $ or
$ F: T_\Omega \to \C^N $. Since component-wise convergence is equivalent to the convergence of mappings, we get that the vector function
$F(z)=(f_1(z),\dots,f_N(z))$ is almost periodic if and only if its components are almost periodic. Therefore for any $\e>0$ the set $E_{\e,K}$
of common almost periods of functions $ f_1,\dots,f_N $ is also  relatively dense.
In particular, this implies that a sum or a product of any finite number of almost periodic functions is also an almost periodic function.

In the rest of the article, only the cases of functions and sets on $\R^d$ or on $T_\Omega \subset C^d$ will be considered.

Now we  give the basic definition of our article.

\begin{Def}
We shall say that functions $f, g$  on $\R^d$   converge weakly at infinity, if
$$
 \lim_{x\to\infty,x\in G} |f(x) - g(x)| = 0,
$$
where $G \subset \R^d$  is a set with the property
\begin{equation}\label{osn}
G\supset\bigcup_{k = 1}^\infty B_\R(x_k ,R_k)\quad\mbox{for some sequence of balls}\quad B_\R(x_k ,R_k),\quad R_k\to\infty.
\end{equation}
\end{Def}
\begin{Def}
 Functions $f, g$ on $T_\Omega$  converge weakly at infinity if for each fixed $y^0\in \Omega$ the functions $f(x + iy^0),\ g(x + iy^0 )$ of the variable $x\in\R^d$ converge weakly at infinity.
\end{Def}
\begin{Th}\label{T3}\hspace{-0,5em}
 If almost periodic functions $f, g$ on $\R^d$ or $T_\Omega$  converge weakly at infinity, then they coincide identically.
\end{Th}

{\bf Proof.} Let $f,g$ be almost periodic functions on $\R^d$.
Fix $x^0  \in \R^d$ and $\e>0$. Let $E_\e$ be the set of $\e$-almost periods of the almost periodic function $h = f - g$.
Taking into account \eqref{osn} and  relative density of $E_\e$, we get that for large $n$  there is a point $\tau_n\in E_\e\cap B_\R(x_n-x^0 ,R_n)$. Hence $x_0+\tau _n\in B_\R(x_n,R_n)$ and $|h(x_0+\tau_n)|<\e$.
Also, $|h(x^0+\tau_n) - h(x^0)|<\e$, therefore, $|h(x^0)|<2\e$. The choice of $\e$ and $x^0$ was arbitrary, therefore $h(x)\equiv0$. In the case of functions on $T_\Omega$ we  take $z^0=x^0 + iy^0\in T_\Omega$ and a compact set $K \subset \Omega$ such that $y^0\in K$, then  replace $E_\e$ by $E_{\e,K}$ and $x^0$ by $z^0$. Theorem is proved.

\medskip

\section{Almost periodic distributions, measures, multisets}\label{S2}

\medskip

Let $D(\R^d)$ be the space of test functions on $\R^d$, i.e., $C^\infty$-functions with compact supports, equipped  with the topology of uniform convergence of derivatives of all orders of functions from $D(\R^d)$, provided that all their supports are subsets of some fixed compact from $R^d$, let $D^*(\R^d)$ be the space of  distributions on $\R^d$, that is, the set of continuous linear functionals on $D(\R^d)$. The distribution space $D^*(T_\Omega)$ is similarly defined as continuous linear functionals on the space $D(T_\Omega)$, consisting of $C^\infty$-functions with compact support in $ T_\Omega$.
\begin{Def}\label{ap}\hspace{-0,5em}
A distribution $f\in D^*(\R^d)$  is called almost periodic, if for any test-function $\p\in D(\R^d)$ the function $(f,\varphi(\cdot-t))$ is almost periodic in the variable $t\in\R^d$.
\end{Def}
\begin{Def}\label{ap1}\hspace{-0,5em}
A distribution $f\in D^*(T_\Omega)$  is called almost periodic, if for any test-function $\p\in D(T_\Omega)$ the function $(f,\p(\cdot-z))$ is almost periodic in the variable $z\in T_\omega$.
Here $\omega$ is the open subset of $\Omega$ such that for all $z\in T_\omega$ the condition $\zeta-z\in\supp\p$ implies $\zeta\in T_\Omega$.
\end{Def}
A particular case of distributions are complex-valued measures. Such measures will be denoted by $\mu$, and the measure, which is the variation of $\mu$, by $|\mu|$. A measure $\mu$ on $\R^d$ is called translation bounded if
$$
\sup_{x\in \R^d} |\mu|(B_\R(x,1)) < \infty.
$$
Similarly, a measure on $T_\Omega$ is called translation bounded if for any compact $K\subset\Omega$
$$
\sup_{x\in \R^d} |\mu|(B_\R(x,1)\times K) \le C,
$$
where the constant $C$ depends on $K$. Note that every nonnegative almost periodic measure is translation bounded. To prove this we should take a nonnegative test function $\varphi(z)\in D(T_\Omega)$ such that $\varphi(z)=1$ on $B_\R(0,1)\times K$,
where $K$ is a compact subset of $\Omega$ (for the case $T_{\R^d}$ we should take nonnegative $\varphi\in D(R^d),\ \varphi(x)=1$ on $B_\R(0,1)$). The function
\begin{equation}
\label{int}
\int {\varphi (z-t)\mu (dz)}
\end{equation}
is almost periodic, hence it is bounded in $t\in\R^d$. On the other hand, for all $t\in\R^d$
$$
\mu(B_\R(t,1)\times K) \le \int\varphi(z-t)\mu(dz).
$$

If a measure $\mu\in D^*(T_\Omega)$  is translation bounded, then we can use any continuous function with compact support as test functions in Definition \ref{ap}.
This follows from the fact that any such a function can be uniformly approximated by $C^\infty$-functions supported on a fixed compact set.
On the other hand, there are signed almost periodic measures for which (\ref{int}) are not almost periodic for an appropriate continuous compactly supported $\varphi$ (\cite{FK2}).
Note that if (\ref{int}) is bounded for all continuous $\varphi$ with compact support, then the complex measure $\mu$ is translation bounded (\cite{R}).

Let $D=\{a,p\}, p\in\N$, be a discrete multiset in $T_\Omega$ or in $\R^d$. It can be identified with a sequence $\{a_n\}$ without condensation points in $T_\Omega$
(or in $\R^d$) such that each point from $T_\Omega$ or in $\R^d$ can occur in this sequence at most a finite number of times. In the case of $ T_\Omega\subset\C$ a discrete multiset is also called
 a divisor (see \cite{FRR}).

\begin{Def} (\cite{FRR})
A discrete multiset $D\subset\R^d$ is called almost periodic if for all $\e>0$ there is a relatively dense set $E_\e=E_\e(D)\subset\R^d$ such that a bijection $\sigma: \N\to\N$ corresponds to any
$\tau\in E_\e$ with the property
$$
\sup_{n\in\N} |a_n-\tau-a_{\sigma (n)}| < \e.
$$
A discrete multiset $D\subset T_\Omega$ is called almost periodic if for all $\e>0$ and compact sets $K\subset\Omega$ there is a relatively dense set $E_{\e,K}=E_{\e,K}(D)\subset\R^d$ such that a bijection
$\sigma: \N\to\N$ corresponds to any $\tau\in E_{\e,K}$ with the property
$$
\sup |a_n-\tau-a_{\sigma (n)}| < \e,
$$
where supremum is taken over all $n\in\N$ such that either $a_n$, or $a_{\sigma(n)}$  belongs to $T_K$.
\end{Def}
We also need a notion of {\it  bounded density}. For a discrete multiset $D\subset\R^d$, $D=\{a_n\}$,  this means that
$$
\sup_{x \in \R^d} \# \{n:\, a_n \in B_\R(x,1) \}< \infty.
$$
 Also, $D\subset T_\Omega$ is of bounded density if for every compact $K\subset \Omega$
\begin{equation}\label{3a}
N(K):=\sup_{x \in \R^d} \#\{n:\, a_n \in B_\R(x,1)\times K\} < \infty.
\end{equation}

It is easy to check that each almost periodic multiset is of bounded density. For $T_{\Omega} \subset \C^d$  the proof can be found in \cite{FRR}.
For convenience, we present it here. The proof for $D \subset \R^d$ differs only in the corresponding simplifications.

Set $\eta= \frac{1}{2} \dist(K, \partial \Omega)$ (in the case $\Omega=\R^d$ set $\eta=\frac{1}{2}$).
Take $R < \infty$ such that every ball $B_\R(x,R)$ intersects with $E_{\eta,K}$. Fix $\tau\in B_\R(x,R)\cap E_{\eta,K}$ and take the bijection $\sigma:\, \N\to\N$ such that for each $a_n\in T_K$
$$
|a_n-\tau-a_{\sigma(n)}|<\eta.
$$
For $a_n \in B_\R(x,1)\times K$ we get
$$
|\Re \ a_{\sigma(n)}|\leq |\Re\ a_{\sigma(n)} - \Re \ a_n + \tau|+ |\Re \ a_n - x| + |x-\tau| < \eta+1+R,
$$
$$
\Im\ a_{\sigma (n)}=\Im\ a_n+(\Im\ a_{\sigma (n)}-\Im\ a_n).
$$
Since $\Im\ a_n\in K$ and $|\Im\ a_{\sigma(n)}-\Im\ a_n|<\eta$, we get $\Im\ a_{\sigma(n)}\in K_1$, where $K_1=\{y:\,\dist(y,K)\le\eta\}$. Thus,
$$
\#\{n: a_n \in B_\R(x,1)\times K\}\le\#\{n: a_\sigma(n)\in B_\R(0,1+\eta+R)\times K_1\},
$$
and we obtain \eqref{3a}.
\medskip

Note that the measure
$$
\mu_D=\sum_n\delta_{a_n},
$$
corresponds to each discrete multiset $D=\{a_n\}$, where $\delta_{a_n}$  is the unit mass at the point $a_n$.

\begin{Th}\label{T4}\hspace{-0,5em}
A discrete multiset $D$ is almost periodic iff the measure $\mu_D$ is almost periodic.
\end{Th}
For $D \subset \C$ this theorem was proved in \cite{FRR}, and for $D\subset~\R^d$ in \cite{FK1}.
Here we give a new, much simpler  proof  for $D \subset T_{\Omega}$. The proof for $D \subset \R^d$ differs only in the corresponding simplifications.
\smallskip

{\bf Proof.}
Let a discrete multiset $D$ be almost periodic and $K\subset\Omega$ be a compact set. Take a function $\varphi\in C^{\infty}(T_{\Omega})$ such that
$\supp \varphi\subset B_\R(0,1/2)\times K$. Let $\e>0$ be arbitrary and  $\delta<(1/2)\dist(K,\partial\Omega)$ such that  for $|z-z'|<\delta$
$$
|\varphi(z)-\varphi(z')|<\frac{\e}{2N(K)},
$$
 where $N(K)$ is defined in \eqref{3a}. Pick $\tau\in E_{\delta,K}(D)$ and the corresponding bijection $\sigma$. We have
$$
\int \varphi(z-\tau)\mu_D(dz)-\int \varphi(z)\mu_D(dz)=\sum_{n}\varphi(a_n-\tau)-\sum_{n}\varphi(a_n)=
$$
$$
=\sum_{n}[\varphi(a_n-\tau)-\varphi(a_{\sigma(n)})].
$$
The number of terms in the letter sum does not exceed $2N(K)$,  moreover, $|a_n-\tau-a_{\sigma(n)}|<\delta$, hence the difference between integrals does not exceed $\e$. Therefore the points of the  set
 $E_{\delta,K}(D)$ are $\e$-almost periods of the function $(\mu_D (\varphi(\cdot-t))$. This reasoning is valid for every $\varphi$ with compact support, therefore the measure $\mu_D$ is almost periodic.

On the other hand, let $\mu_D$ be the almost periodic measure on $T_{\Omega}$, which corresponds to a discrete multiset $D=\{a_n\}$.  Fix a compact set $K\subset \Omega$ and $\e<\frac{1}{4}\min \{1,\dist(K,\Omega)\}$. Put
$$
\widetilde{K}=\{y\in \Omega : \dist(y,K)\leq\e\}.
$$
Choosing a sufficiently large $K$, we can assume that either $D\subset K$, or $D\setminus T_{\tilde K}\neq\emptyset$. Since $\mu_D$ is almost periodic we get that it is translation bounded, hence for some $N<\infty$
$$
\mu_D(B_\R(x,1)\times \widetilde{K})\leq N,\quad\forall x\in \R^d,
$$
therefore,
\begin{equation}
\label{N}
\# \{n:a_n\in B_\R(x,1)\times \widetilde{K}\}\le N, \,\forall x\in \R^d.
\end{equation}
Set $\delta=\e/(4N+1)$. Let $A$ be any connected component of the set $\bigcup\limits_{n} B_\C(a_n, 2\delta)$ such that $A\cap T_K\neq\emptyset$. There exists $a_{n'}\in A$ such that $ B_\C(a_{n'}, 2\delta)\cap T_K\neq\emptyset$. If  $A\cap\partial B_\C(a_{n'},\e)\neq\emptyset$, then the connected set $A\cap B_\C(a_{n'}, \e)$ contains at least  $\e/(4\delta)>N$ points of $D$, which contradicts (\ref{N}).
Hence,
$$
A\subset B_\C(a_{n'},\e)\subset B_\R(\Re\,a_{n'},1)\times \widetilde{K}
$$
and, by (\ref{N}), $\#\{n:a_n\in A\}\le N$.

By $\varphi(z)$ denote any  $C^{\infty}$-function on  $\C^d$  such that
$$
0\le \varphi(z)\le 1,\quad \varphi(0)=1, \quad \supp\varphi\subset B_\C(0,1),
$$
Let $\a=\int\p(z)\omega(dz)$, where $\omega$ is the Lebesgue measure on $\C^d$.  Put
$$
\Psi(z):=\int\varphi\left(\frac{z-w}{\d}\right)\mu_D(dw)=\sum_n\varphi\left(\frac{z-a_n}{\d}\right).
$$
Since
$$
\dist({\widetilde{K}},\partial\Omega)\ge \dist(K,\partial\Omega)-\e\ge \e,
$$
 we see that $\Psi(z)$ is defined and almost periodic on $T_{\widetilde{K}}$.
Let $\tau$ be $\rho$-almost period of $\Psi(z)$ with $\rho<\min\{1;2^{-2d}\a/(N\omega_{2d})\}$, where $\omega_{2d}=\omega(B_\C(0,1))$.   We have
\begin{equation}
\label{fe}
|\Psi(z+\tau)-\Psi(z)|<\rho, \quad \forall z\in T_{\widetilde{K}}.
\end{equation}
On the other hand,
$$
\Psi(z)=0\quad\mbox{for}\quad z\notin\cup_n B_\C(a_n,\delta)\quad\mbox{and}\quad\Psi(z+\tau)=\Psi(a_n)\ge1\quad\mbox{for}\ z=a_n-\tau.
$$
Therefore the set $A\setminus\cup_n B_\C(a_n,\delta)$ does not contain any point $a_n-\tau$. If $A'$ is
another connected component of the set $\cup_n B_\C(a_n,2\delta)$,
then for the same reason $A'\setminus\cup_n B_\C(a_n,\delta)$ does not contain any point $a_n-\tau$.
 Thus the set $A$ contains all  balls  $B_\C(a_n,\delta)$, for which $a_n\in A$ and all balls $B_\C(a_n-\tau,\delta)$, for which $a_n-\tau\in A$, and does not intersects with balls $B_\C(a_n,\delta)$ with
$a_n\notin A$ and balls $B_\C(a_n-\tau,\delta)$ with $a_n-\tau\notin A$.
We get
$$
\a\d^{2d}\#\{n : a_n\in A\}=\sum_{n:a_n\in A}\int\varphi\left(\frac{z-a_n}{\d}\right)\omega(dz)=\int_A\Psi(z)\omega(dz),
$$
$$
\alpha\d^{2d}\#\{n : a_n-\tau\in A\}=\sum_{n:a_n-\tau\in A}\int\varphi\left(\frac{z+\tau-a_n}{\d}\right)\omega(dz)=\int_A\Psi(z+\tau)\omega(dz),
$$
Note that
$$
\int_A \omega(dz)\le\sum_{a_n\in A}\int_{B_\C(a_n,2\d)}\omega(dz)= N\omega_{2d}(2\d)^{2d}.
$$
By \eqref{fe},
$$
|\#\{n:a_n\in A\}-\#\{n:a_n-\tau\in A\}|\le\frac{\int_A|\Psi(z)-\Psi(z+\tau)|\omega(dz)}{\d^{2d}\a}<\frac{\rho N\omega_{2d}2^{2d}}{\a}<1.
$$
Therefore,
$$
\#\{n:a_n\in A\}=\#\{n:a_n-\tau\in A\},
$$
 which allows to construct a bijection $\sigma$ between the sets
$\{n:a_n-\tau\in A\}$ and $\{n~:~a_n\in~A\}$.

This construction works for every connected component of the set $\cup_n B_\C(a_n, 2\delta)$, hence there exists a bijection $\sigma$ of a part $S_1$ of $\N$ to the part $S_2$ of $\N$.
It follows from the inequality  $\diam A\leq2\e$ that
\begin{equation}\label{fn}
|a_n-\tau-a_{\sigma(n)}|<2\e.
\end{equation}
If $D\subset T_K$, we have $S_1=S_2=\N$, and theorem is proved. If $D\setminus T_{\widetilde{K}}\neq\emptyset$, we have only
$$
\{n:a_n\in T_{K}\}\subset S_1\cup S_2\subset\{n:~a_n\in T_{\widetilde{K}}\}.
$$
 For $a\in D\setminus T_{\tilde K}$ put
 $$
 \eta<\frac 12 \min\{\dist (\Im\ a,\widetilde{K}), \dist \{\Im\ a,\partial\Omega\}\}
 $$
 and consider the function
 $$
\Psi(z)=\int\varphi\left(\frac{w-z}{\eta}\right)\mu_D(dw)=\sum_{n:a_n\in D}\varphi\left(\frac{a_n-z}{\eta}\right).
$$
In view of the choice of $\eta$, this function is well-defined and almost periodic on $T_\omega$ with $\omega=\{y\in\Omega:\, \dist(y,\partial\Omega)>\eta\}$.
Furthermore,  $\Psi(a)\ge \varphi(0) = 1$, hence  $\Psi(a+t)$ is strictly positive for some  large $t\in\R^d$. Therefore the set $\{n:a_n\in D\backslash\widetilde{K}\}$ is unbounded and countable,
as well as the sets $\N\setminus S_1$ and $\N\setminus S_2$. For points $a_n$ with $n\notin S_1$ the condition (\ref{fn}) need not be required, therefore the bijection $\sigma : S_1\to S_2$
can be extended to a bijection $\N\to\N$. The Theorem is proved.

\medskip

\section{Uniqueness theorems for almost periodic distributions, measures, multisets}\label{S3}

\medskip
\begin{Def}
We shall say that distributions $f, g \in D^*(\R^d)$  converge weakly at infinity, if for any $\varphi\in D(\R^d)$
the functions $(f,~\varphi~(\cdot~-t))$ and $(g,~\varphi(~\cdot-t))$ of the variable $t\in \R^d$  converge weakly at infinity.

Also, we shall say that distributions $f, g \in D^*(T_{\Omega})$ converge weakly at infinity, if for any $\p\in D(T_{\Omega})$ the functions $(f,\p(\cdot-z))$ and $(g,\p(\cdot-z))$
of the variable $z\in T_\omega$ converge weakly at infinity ($\omega\subset\Omega$ is defined in Definition \ref{ap1}).
\end{Def}

It follows from Theorem \ref{T3}
\begin{Th}\label{T5}
If two almost periodic distributions or measures $f, g \in D^*(\R^d)$  converge weakly at infinity, then $f\equiv g$.
The similar assertion is valid for $f, g \in D^*(T_{\Omega})$.
\end{Th}

\begin{Def}
We shall say that two discrete multisets $F=\{a_n\}, H=\{b_n\} \subset \R^d$  converge weakly at infinity, if there is a set
$G\subset \R^d$ satisfying (\ref{osn}) such that under an appropriate numbering
$$
\lim_{n\rightarrow\infty, n\in N(G)}{a_n-b_n}=0,
$$
где $\N(G)=\{n\in\N : a_n\ \mbox{ or }\ b_n\in G \}$.
\end{Def}

\begin{Def}
\label{def7}
We shall say that two discrete multisets $F=\{a_n\}, H=\{b_n\} \subset T_\Omega$  converge weakly at infinity, if  for every
$K \subset \Omega$ there is a set $G=G(K) \subset \R^d$ satisfying (\ref{osn}) such that under an appropriate numbering
$$
\lim_{n\rightarrow\infty, n\in\N(G,K)}{a_n-b_n}=0,
$$
where $\N(G,K)=\{n\in\N : a_n\in G \times K\ \mbox{ or }\ b_n\in G\times K\}$.
\end{Def}
\begin{Th}\label{T6}
If two discrete multisets $F=\{a_n\}, H=\{b_n\}$  converge weakly at infinity, then they are identical.
\end{Th}

{\bf Proof}.
It follows from theorems \ref{T4} and \ref{T5} that we have to check the weak convergence of measures  $\mu_F$ and $\mu_H$ at infinity. The latter means that for any $\varphi\in D(\R^d)$ (or $\varphi~\in~D(T_{\Omega})$) the almost periodic functions of the variable $t\in\R^d$
$$
\Psi_F(t)=(\mu_F, \varphi(\cdot-t))=\sum\limits_{n}\varphi(a_n-t)
$$
and
$$
\Psi_H(t)=(\mu_H,\varphi(\cdot-t))=\sum\limits_{n}\varphi(b_n-t)
$$
 converge weakly at infinity. To be specific consider the case $F, H \subset T_{\Omega}$. A proof is the same for $F,H \subset \R^d$.

Suppose that $\supp \varphi\subset B_\R(0,1)\times K$ for compact $K\subset\Omega$.
Take $\e>0$ and then $\delta >0$ such that $|\varphi(z)-\varphi(z')|<\e/(N(K))$ for $|z-z'|<\delta$, where $N(K)$ is the constant from (\ref{3a}). Let a set $G\subset\R^d$  satisfy \eqref{osn}
with balls $B_\R(x_k,R_k),\ k\in\N$. It is easy to see that having reduced by 3 times the radii of these balls and changing the location of their centers, we can assume that $\dist(B(x_k,R_k),0)\to\infty$.
 For sufficiently large $k$ and for $a_n,\ b_n\in B_\R(x_k,R_k)\times K$ we have $|a_n-b_n|<\delta$. Also assume that $R_k>2$.

Let $t\in B_\R(x_k, R_k/2)$ and $a_n-t\in \supp \varphi$. Then $a_n \in B_\R(x_k,R_k)\times K$, and the same is valid for $b_n-t$.
Therefore if  $a_n-t\in \supp \varphi$ or $b_n-t\in \supp \varphi$, we get $|a_n-b_n|<\delta$ and
$$
|\Psi_F(t)-\Psi_H(t)|\le \sum_n|\varphi(a_n-t)-\varphi(b_n-t)|<\frac {\e}{N(K)}\cdot N(K)=\e.
$$
Hence the almost periodic functions $\Psi_F(t)$ и $\Psi_H(t)$  converge weakly at infinity.

\medskip

\section{Uniqueness theorems for delta-subharmonic and meromorphic functions}\label{S4}

\medskip

It follows immediately from the definition that any partial derivative of an almost periodic distribution from
$D^*(\R^d)$ or $D^*(T_{\Omega})$ is also an almost periodic distribution. Since any subharmonic function on any region from $\R^d$ is locally integrable, it can be considered as a distribution. Thus, if $u$ is a subharmonic almost periodic function on $D^*(\R^d)$ or $D^*(T_{\Omega})$,  then its Riesz measure $\Delta u$ is also an almost periodic distribution, and the same is true for the difference of subharmonic functions, the so-called delta-subharmonic functions.

It follows from Theorem \ref{T5}

\begin{Th}\label{T7}
If two delta-subharmonic functions $u, v$ on $\R^d$ or $T_{\Omega}$ have  converging weakly at infinity Riesz measures $\Delta u$ and $\Delta v$, then $u=v+h$ with a harmonic function $h$.
\end{Th}

The last part of the proof uses the fact that the condition $\Delta h=0$ in the sense of distributions implies that $h$ is an ordinary harmonic function.

\begin{Def}(см. \cite{S}, \cite{FP})
A meromorphic function $f(z)$ on the strip $S_{a,b}=\{z\in\C: \Re\ z\in\R, a<\Im\ z<b\}$,  $-\infty \le a<b\le +\infty$, is called almost periodic, if in any smaller strip $S_{\alpha,\beta}$, $a<\alpha<\beta<b$, the function $\rho_S(f(z+t),f(z))$, where $\rho_S$ is the spherical distance, is almost periodic in the variable $t\in\R$.
\end{Def}

In \cite{FP} the following properties of meromorphic almost periodic functions are proved:
\begin{itemize}

 \item The distance between any pole and any zero of meromorphic almost periodic functions is bounded from below by a strictly positive constant depending on the strip in which this pole and zero lie,

 \item Every meromorphic almost periodic function on $S_{a,b}$ is a ratio of two holomorphic almost periodic functions in $S_{a,b}$; the converse assertion is only valid if distances between poles and zeros of this ratio are uniformly bounded from below by a strictly positive constant in any smaller strip. In particular, every holomorphic almost periodic function in a strip is simultaneously a meromorphic almost periodic function.

\end{itemize}

 It was proved in \cite{R} that for any almost periodic holomorphic function $f$ on $S_{a,b}$ the function $\log|f|$ is an almost periodic distribution, hence the measure $\mu_Z$
 corresponding to the multiset of zeros $Z_f$ is almost periodic.  Also, if $f$ is an  almost periodic meromorphic function, then the measures $\mu_Z$ and $\mu_P$ corresponding to the multiset of zeros $Z_f$ and the multiset of poles $P_f$ of $f$  are also almost periodic. Therefore, Theorem \ref{T7} implies

\begin{Th}\label{T8}
If multisets of poles $P_f$ and $P_g$ of meromorphic almost periodic functions $f, g$ in a strip $S_{a,b}$  converge weakly and the same is true for multisets of zeros
 $Z_f$ and $Z_g$, then $P_f=P_g$, $Z_f=Z_g$, hence, $f/g$ is a holomorphic almost periodic function on $S_{a,b}$ without zeros.
\end{Th}

If $f, g$ are holomorphic almost periodic functions on $S_{a,b}$ and multisets of zeros
 $Z_f,\ Z_g$  converge weakly at infinity, then $Z_f=Z_g$, and we obtain Theorem 6 from \cite{KS}.

Note that a linear-fractional mapping of a meromorphic almost periodic function $f$ is a meromorphic almost periodic function too. Then instead of zeros and poles one can consider $A_1$ -points and
$A_2$ -points, $A_1\neq A_2$, that is zeros of functions $f-A_1$ and $f-A_2$.
Also, for $T_\Omega=\C$ we obtain the following theorem:

\begin{Th}\label{T9}\hspace{-0,5em}
Let $f, g$ be meromorphic almost periodic functions on $\C$ and let $A_j$-points of $f$ converge weakly at infinity  to $A_j$-points of $g$ for  three pairwise distinct values
$A_1,\,A_2,\,A_3$. Then either $f=g$, or $f$ and $g$ have the forms
\begin{equation}\label{e}
  f=T\left(\frac{1-h_1}{h_2-h_1}\right),\qquad g=T\left(\frac{h_2-h_1h_2}{h_2-h_1}\right),
\end{equation}
where $h_1,\,h_2$ are distinct entire functions without zeros, and $T$ is a linear-fractional mapping that moves the triple point $0,1,\infty$ to the triple point $A_1,\,A_2,\,A_3$.
\end{Th}
At the final stage of the proof we use the following theorem from \cite{N}:
If two meromorphic functions on $\C$ have the same multisets of $A$ -points for three distinct values of $A_1,\,A_2,\,A_3$, then these functions either coincide,
 or have  form \eqref{e}.

\end{document}